# FAMILIES OF CURVES AND WEIGHT DISTRIBUTIONS OF CODES

RENÉ SCHOOF

ABSTRACT. In this expository paper we show how one can, in a uniform way, calculate the weight distributions of some well-known binary cyclic codes. The codes are related to certain families of curves, and the weight distributions are related to the distribution of the number of rational points on the curves.

## 1. INTRODUCTION

Let $\mathbf{F}_q$ be a finite field with $q$ elements. A *code over* $\mathbf{F}_q$ is simply a linear subspace
$$C \subset \mathbf{F}_q^n.$$
The most natural and interesting codes are *binary,* i.e. they are codes over $\mathbf{F}_2$. These codes are the principal object of study in this paper, but we will also encounter codes over extension fields of $\mathbf{F}_2$.

Around 1980 the Russian mathematician and engineer V.G. Goppa discovered a beautiful connection between codes and algebraic curves over finite fields. He obtained a very general method to construct codes from certain linear systems on algebraic curves.

In this expository paper we discuss a *different* relation between codes and algebraic curves. We study certain families of binary codes whose words are parametrized by algebraic curves. These codes are very classical, and many of their properties are derived in the famous book by MacWilliams and Sloane [8]. There the methods are of a combinatorial nature. Our approach is entirely geometric. We encounter elliptic curves and Hecke operators for the modular group $\Gamma_1(4)$. Our intended reader has some familiarity with the theory of elliptic curves and is not supposed to know *anything* about coding theory, so we will explain some very well-known coding theoretical concepts and say nothing about some deep results concerning elliptic curves over finite fields.

The elements of a code are called *code words*; the dimension $n$ is called the *length* of the code. The *weight* $|\mathbf{v}|$ of a code word $\mathbf{v} = (v_1, v_2, \ldots, v_n)$ is defined by
$$|\mathbf{v}| = \#\{i : v_i \neq 0\}.$$
Coding theorists are usually interested in codes which are "good", i.e. for which the *minimum distance*
$$d = \min\{|\mathbf{v}| : \mathbf{v} \in C, \ \mathbf{v} \neq 0\}$$
is large with respect to $\dim(C)$. In this paper we do not take this point of view. We merely investigate the distribution of the weights in a given code.

An easy example of a binary code is the *repetition code*, which consists of the zero vector and the vector $(1, 1, \ldots, 1)$. It has length $n$, dimension 1, and minimum distance $d = n$. Another example is the *parity check code*,

$$\mathrm{P} = \{\mathbf{v} \in \mathbf{F}_2^n : \sum_{i=0}^{n} v_i = 0\}.$$

It has length $n$, dimension $n-1$, and minimum distance $d = 2$.

An invariant of a code $C$ which is finer than the minimum distance is its *weight distribution*. It is the function

$$A_i = \#\{\mathbf{v} \in C : |\mathbf{v}| = i\}, \qquad 0 \leq i \leq n,$$

and it is often given in terms of a polynomial, the *weight enumerator* $W_C(X) \in \mathbf{Z}[X]$:

$$W_C(X) = \sum_{i=0}^{n} A_i X^i,$$
$$= \sum_{\mathbf{v} \in C} X^{|\mathbf{v}|}.$$

The minimum distance is the smallest integer $d > 0$ for which $A_d \neq 0$. The weight enumerator is the finite analogue of the $\vartheta$-series associated to an integral lattice in Euclidean space.

It is clear that all these notions depend in an essential way on the basis of the vector space $\mathbf{F}_2^n$. Changing the basis changes everything. Only permuting the basis vectors does not change the code "in an essential way".

Determining the weight distribution of a code is, in general, a difficult task. In this paper we discuss four families of codes for which, luckily, one is able to determine the weight distributions. These codes all belong to a classical, well-studied class of codes: they are all *cyclic*.

**Definition.** A cyclic code $I$ is an ideal

$$I \subset \mathbf{F}_2[T]/(T^n - 1),$$

where the vector space $\mathbf{F}_2[T]/(T^n - 1)$ is equipped with the $\mathbf{F}_2$-basis consisting of the monomials $1, T, \ldots, T^{n-1}$. In other words, the coordinates $(v_0, v_1, \ldots, v_{n-1})$ of the code words in $I$ are simply the coefficients of the polynomials $v_0 + v_1 T + \cdots + v_{n-1} T^{n-1} \in \mathbf{F}_2[T]/(T^n - 1)$.

The codes $I$ are called cyclic since the word $(v_1, \ldots, v_{n-1}, v_0)$ is contained in $I$ whenever $(v_0, v_1, \ldots, v_{n-1}) \in I$.

For instance, taking for $I$ the ideal generated by the polynomial $T - 1$, we find

$$I = \{v_0 + v_1 T + \cdots + v_{n-1} T^{n-1} \in \mathbf{F}_2[T]/(T^n - 1) : v_0 + v_1 + \cdots + v_{n-1} = 0\}.$$

This is precisely the parity check code P. One verifies easily that the ideal generated by the polynomial $T^{n-1} + \cdots + T + 1$ is simply the repetition code.

Another example is the family of *Hamming codes H:* let $n = q - 1$ where $q = 2^m$ for some $m \geq 2$. Then

$$\mathrm{H}_m \subset \mathbf{F}_2[T]/(T^{q-1} - 1)$$

is the ideal generated by the minimum polynomial $f_\alpha(T) \in \mathbf{F}_2[T]$ of $\alpha \in \mathbf{F}_q^*$ where $\alpha$ is a generator of the multiplicative group $\mathbf{F}_q^*$. We will see later that, in a precise sense, the codes do not depend on the choice of $\alpha$.

In this paper we consider for $m > 2$ the so-called *double error correcting BCH code*
$$\mathrm{B}_m \subset \mathbf{F}_2[T]/(T^{q-1} - 1);$$
this is the ideal generated by the product $f_\alpha(T) f_{\alpha^3}(T)$. Finally we consider the *Melas codes*
$$\mathrm{M}_m \subset \mathbf{F}_2[T]/(T^{q-1} - 1),$$
which are the ideals generated by the product $f_\alpha(T) f_{\alpha^{-1}}(T)$.

It is, of course, perfectly possible to study weight distributions without ever speaking about codes. The weight distribution of a cyclic code $I$ simply carries information about the number of polynomials in the ideal $I$ with a given number of non-zero coefficients. For instance, for a cyclic code $I$ of length $n$ one has that
$$A_i = \#\{f \in I : \deg(f) < n \text{ and } f \text{ has precisely } i \text{ non-vanishing coefficients}\}.$$

In this paper we determine the weight distributions of the codes $\mathfrak{P}$, $\mathfrak{H}_m$, $\mathfrak{B}_m$, and $\mathfrak{M}_m$. None of this is new. For the codes $\mathfrak{P}$, $\mathfrak{H}_m$, and $\mathfrak{B}_m$ the results can be found in the book by MacWilliams and Sloane [8]; the Melas codes $\mathfrak{M}_m$ are dealt with in [12]. The algebraic geometric approach, however, is different from the discussion in [8]. We relate the codes $\mathfrak{H}_m$, $\mathfrak{B}_m$, and $\mathfrak{M}_m$ to certain families of curves and their weight distributions to the distribution of the number of rational points on the curves in the family.

In section 2 we explain the *MacWilliams Identities,* relating the weight distribution of a code to the one of its dual. As an application we determine the weight distribution of the parity check codes $\mathfrak{P}$. In section 3 we determine the weight distribution of the Hamning codes $\mathfrak{H}_m$. This is done by extending the code to $\mathbf{F}_q$ and by applying *Delsarte's Theorem.* In section 4 we study the BCH codes $\mathfrak{B}_m$. This involves a family of supersingular elliptic curves over $\mathbf{F}_q$. We determine how often every isomorphism class occurs in the family, and this gives us the weight distribution of the codes $\mathfrak{B}_m$. The calculation is a special case of the families of curves studied in [14]. Finally in section 5 we determine the weight distribution of the Melas codes $\mathfrak{M}_m$. This time a family of non-supersingular elliptic curves over $\mathbf{F}_q$ is involved. The weight distributions of the dual codes can be described in terms of certain class numbers of binary quadratic forms. An application of the Eichler-Selberg Trace Formula gives, in a somewhat mysterious way, the weight distribution of $\mathfrak{M}_m$. The result involves the traces of the Hecke operators $T_q$ acting on the cusp forms for the modular group $\Gamma_1(4)$. This is probably related to the fact that the curves in our family all admit a rational point of order 4, but we do not establish a direct connection.

It would be very interesting to extend the methods of this paper to other families of cyclic codes. This seems difficult since it involves, in general, curves of genus larger than 1 and difficult questions about their moduli and their numbers of rational points. See [14] where Reed-Muller codes are discussed from this point of view.

## 2. The parity check codes

Often the determination of the weight distribution of a code $C$ goes together with the determination of the weight distribution of the *dual code* $C^\perp$. For instance, if $C = C^\perp$, i.e. if $C$ is self dual, much can be said about the weight distribution of $C$. In this case there are beautiful connections with invariant theory and the theory of $\vartheta$-functions. See [8, Ch. 19].

**Definition.** Let $C \subset \mathbf{F}_q^n$ be a code. The *dual code* $C^\perp$ is defined by
$$C^\perp = \{\mathbf{v} \in \mathbf{F}_q^n : \langle \mathbf{v}, \mathbf{w} \rangle = 0 \text{ for all } \mathbf{w} \in C\}$$
where for $\mathbf{v} = (v_1, v_2, \ldots, v_n)$ and $\mathbf{w} = (w_1, w_2, \ldots, w_n)$ in $\mathbf{F}_q^n$ the scalar product $\langle \mathbf{v}, \mathbf{w} \rangle$ is defined by
$$\langle \mathbf{v}, \mathbf{w} \rangle = \sum_{i=1}^n v_i w_i.$$

It is easily seen that the dual of $C^\perp$ is the code $C$ again. For example, the parity check code P and the repetition code are dual to one another.

The following theorem gives, for binary codes, a relation between the weight distribution of the codes $C$ and $C^\perp$.

**Theorem 2.2** (The MacWilliams Identities). *Let $C$ be a binary code of length $n$ with weight distribution $W_C(X)$. Then*
$$W_{C^\perp}(X) = \frac{1}{\#C} \sum_{i=0}^n A_i (1-X)^i (1+X)^{n-i}.$$

*Proof.* For a function $f$ on $\mathbf{F}_2^n$ we define its Fourier transform $\hat{f}$ by
$$\hat{f}(\mathbf{v}) = \sum_{\mathbf{w} \in \mathbf{F}_2^n} (-1)^{\langle \mathbf{v}, \mathbf{w} \rangle} f(\mathbf{w}), \quad \text{for } \mathbf{v} \in \mathbf{F}_2^n.$$

Since $\sum_{\mathbf{v} \in C} (-1)^{\langle \mathbf{v}, \mathbf{w} \rangle}$ is equal to $\#C$ or $0$ depending on whether $\mathbf{w} \in C^\perp$ or not, one deduces the following inversion formula
$$\sum_{\mathbf{v} \in C} \hat{f}(\mathbf{v}) = \#C \sum_{\mathbf{w} \in C^\perp} f(\mathbf{w}).$$

We apply this to the function
$$f(\mathbf{w}) = X^{|\mathbf{w}|}.$$

One has that
$$W_{C^\perp}(X) = \sum_{\mathbf{w} \in C^\perp} X^{|\mathbf{w}|} = \frac{1}{\#C} \sum_{\mathbf{v} \in C} \hat{f}(\mathbf{v}).$$

Since
$$\hat{f}(\mathbf{v}) = \sum_{\mathbf{w} \in \mathbf{F}_2^n} (-1)^{\langle \mathbf{v}, \mathbf{w} \rangle} X^{|\mathbf{w}|} = (1-X)^{|\mathbf{v}|}(1+X)^{n-|\mathbf{v}|},$$

the theorem follows easily. The last fact is proved by induction with respect to $n$. This completes the proof of the theorem.

As an easy example, we determine the weight distribution of the parity check codes P. These codes are dual to the repetition codes which are easily seen to have the following weight enumerator
$$W(X) = 1 + X^n.$$

By Theorem 2.1 we have therefore that
$$W_\mathfrak{P}(X) = \frac{1}{2}\left((1-X)^n + (1+X)^n\right),$$
$$= \sum_{\substack{i=0 \\ i \text{ odd}}}^n \binom{n}{i} X^i,$$

a result which is also easily obtained directly.

## 3. The Hamming codes

Let $m \geq 2$ and $q = 2^m$. The Hamming code $\mathfrak{H}_m$ of length $q-1$ is the ideal in the ring $\mathbf{F}_2[T]/(T^{q-1}-1)$ generated by the polynomial $f_\alpha(T)$. Here $f_\alpha(T) \in \mathbf{F}_2[T]$ denotes the minimum polynomial of a generator $\alpha$ of the multiplicative group $\mathbf{F}_q^*$. As we explained in the introduction, we take the monomials $1, T, \ldots, T^{q-2}$ as a basis for the vector space $\mathbf{F}_2[T]/(T^{q-1}-1)$. We will see below that the choice of $\alpha$ is not important.

**Proposition 3.1.** *The dimension of $\mathfrak{H}_m$ is $q-1-m$ and its minimum distance is $3$.*

*Proof.* Since $\mathbf{F}_2[T]/(f_\alpha(T))$ is a field of $2^m$ elements, the ideal $\mathfrak{H}_m = (f_\alpha(T))$ itself has $2^{q-1}/2^m$ elements. Therefore the dimension of $\mathfrak{H}_m$ is $q-1-m$.

Suppose $h \in \mathfrak{H}_m$ is a non-zero word of weight at most 2. In other words, $h(T) = T^a$ for some $0 \leq a < q-1$ or $h(T) = T^a + T^b$ for some $0 \leq a < b < q-1$. Since $h(\alpha) = 0$, the first possibility cannot occur. In the second case we have $\alpha^{b-a} = 1$ which is impossible since $\alpha$ has order $q-1$. Therefore the minimum distance is at least 3.

This implies that the "unit balls"
$$S(\mathbf{v}) = \{\mathbf{w} \in \mathbf{F}_2^{q-1} : |\mathbf{v} - \mathbf{w}| \leq 1\}$$
with "center" $\mathbf{v} \in \mathfrak{H}_m$ are mutually disjoint. Since each ball $S(\mathbf{v})$ has cardinality $2^m$, their union has cardinality $\#\mathfrak{H}_m 2^m = 2^{q-1-m} 2^m = 2^{q-1}$. This is precisely the cardinality of the vector space $\mathbf{F}_2[T]/(T^{q-1}-1)$, and we conclude that the balls $S(\mathbf{v})$ cover the vector space $\mathbf{F}_2^{q-1}$. This easily implies that the minimum distance is actually *equal* to 3: the polynomial $1 + T$ is contained in a unit ball $S(\mathbf{v})$ for some code word $\mathbf{v}$. Therefore $1 + T - \mathbf{v}$ is a monomial $T^a$ for some $a \neq 0, 1$, and we see that $T^a + T + 1 = \mathbf{v}$ is a code word of weight 3. This proves the proposition.

Taking $m = 3$ and $f_\alpha(T) = T^3 + T + 1$, one obtains the famous Hamming code of length 7 and dimension 4. Here is a complete list of all 16 vectors, with respect to the basis $1, T, \ldots, T^6$:

$$\begin{array}{cccccccccccccccc}
0 & 0 & 0 & 0 & 1 & 1 & 0 & 1 & 1 & 1 & 1 & 1 & 0 & 0 & 1 & 0 \\
0 & 0 & 0 & 1 & 1 & 0 & 1 & 0 & 1 & 1 & 1 & 0 & 0 & 1 & 0 & 1 \\
0 & 0 & 1 & 1 & 0 & 1 & 0 & 0 & 1 & 1 & 0 & 0 & 1 & 0 & 1 & 1 \\
0 & 1 & 1 & 0 & 1 & 0 & 0 & 0 & 1 & 0 & 0 & 1 & 0 & 1 & 1 & 1 \\
0 & 1 & 0 & 1 & 0 & 0 & 0 & 1 & 1 & 0 & 1 & 0 & 1 & 1 & 1 & 0 \\
0 & 0 & 1 & 0 & 0 & 0 & 1 & 1 & 1 & 1 & 0 & 1 & 1 & 1 & 0 & 0 \\
0 & 1 & 0 & 0 & 0 & 1 & 1 & 0 & 1 & 0 & 1 & 1 & 1 & 0 & 0 & 1 \\
\end{array}$$

To study the codes $\mathfrak{H}_m$, we extend the base field to $\mathbf{F}_q$:
$$\mathbf{F}_2[T]/(T^{q-1}-1) \subset \mathbf{F}_q[T]/(T^{q-1}-1).$$

Since the words in $(T-\alpha)$ are the polynomials $a_{q-2}T^{q-2} + \cdots + a_1 T + a_0$ that vanish at $\alpha$, we see that $\mathfrak{H}_m$ is the *restriction* of the ideal $(T-\alpha) \subset \mathbf{F}_q[T]/(T^{q-1}-1)$, i.e.
$$\mathfrak{H}_m = (T-\alpha) \cap \mathbf{F}_2[T]/(T^{q-1}-1).$$

It is easy to compute the dual of the $\mathbf{F}_q$-code $(T-\alpha)$. The polynomials in $(T-\alpha)$ are the vectors $(a_0, a_1, \ldots, a_{q-2})$ that are orthogonal to $(1, \alpha, \alpha^2, \ldots, \alpha^{q-2})$.

Since $\alpha$ generates $\mathbf{F}_q^*$, the latter vector has precisely all elements $x \in \mathbf{F}_q^*$ as its coordinates. Therefore the $\mathbf{F}_q$-dual of $(T - \alpha)$ is given by

$$(T - \alpha)^\perp = \{(\lambda x)_{x \in \mathbf{F}_q^*} : \lambda \in \mathbf{F}_q\}.$$

We note in passing that this dual code does not depend on $\alpha$. Therefore, up to a permutation of the basis, the Hamming code $\mathfrak{H}_m$ does not depend on the choice of $\alpha$ either, but only on $m$.

We wish to compute the weight distributions of the codes $\mathfrak{H}_m$, which themselves are the restrictions of the codes $(T - \alpha)$. By the MacWilliams Identities, the following theorem reduces this to a study of the *Trace code* of $(T - \alpha)^\perp$:

**Theorem 3.2** (Delsarte's Theorem). *Let $\mathbf{F}_q$ be an extension of $\mathbf{F}_2$, and let $C$ be a code over $\mathbf{F}_q$. Then*

$$\mathrm{Trace}(C)^\perp = \mathrm{Res}(C^\perp),$$

*where $\mathrm{Res}(C^\perp)$ is the restriction of $C^\perp$ to $\mathbf{F}_2$ and $\mathrm{Trace}(C)$ is the $\mathbf{F}_2$-code given by*

$$\mathrm{Trace}(C) = \{\mathop{\mathrm{Tr}}_{\mathbf{F}_q/\mathbf{F}_2}(\mathbf{v}) : \mathbf{v} \in C\}.$$

*Proof.* By definiton $\mathbf{v} \in \mathrm{Trace}(C)^\perp$ if and only if $\langle \mathbf{v}, \mathrm{Tr}(\mathbf{w}) \rangle = 0$ for all $\mathbf{w} \in C$. Since $\mathbf{v}$ has coordinates in $\mathbf{F}_2$, this is equivalent to

$$\mathrm{Tr}(\langle \mathbf{v}, \mathbf{w} \rangle) = 0 \quad \text{for all } \mathbf{w} \in C.$$

This is also true for $\lambda \mathbf{w}$ with $\lambda \in \mathbf{F}_q$ instead of $\mathbf{w}$. Since the trace is non-degenerate, we conclude therefore that the last statement is equivalent to $\langle \mathbf{v}, \mathbf{w} \rangle = 0$ for all $\mathbf{w} \in C$, i.e. to $\mathbf{v} \in C^\perp$ and the proof is complete.

So, we have the following situation:

$$\begin{array}{ccc} \mathfrak{H}_m & \overset{\mathrm{Res}}{\longleftarrow} & (T - \alpha) \\ \Big\updownarrow \mathrm{dual} & & \Big\updownarrow \mathrm{dual} \\ \mathfrak{H}_m^\perp & \overset{\mathrm{Trace}}{\longleftarrow} & (T - \alpha)^\perp \end{array}$$

Because of the simple description of $(T - \alpha)$ and $(T - \alpha)^\perp$, we can compute the duals of the Hamming codes: by the discussion above and Delsarte's Theorem we have

$$\mathfrak{H}_m^\perp = \{\mathrm{Tr}(\lambda x)_{x \in \mathbf{F}_q^*} : \lambda \in \mathbf{F}_q\}.$$

The cardinality of $\mathfrak{H}_m^\perp$ is $q$. It is very easy to determine the weight distribution of $\mathfrak{H}_m^\perp$: for $\lambda = 0$ we get the zero word of weight 0 and for $\lambda \neq 0$ we get a word of weight $q/2$, because $\lambda x$ runs through all non-zero elements of $\mathbf{F}_q$.

By the MacWilliams Identities we have that

$$W_{\mathfrak{H}_m}(X) = \frac{1}{q}\left((1+X)^{q-1} + (q-1)(1-X)^{q/2}(1+X)^{q/2-1}\right);$$

explicitly, the number of code words in $\mathfrak{H}_m$ of weight $i$ is given by

$$A_i = \frac{1}{q}\left(\binom{q-1}{i} + (q-1)(-1)^{[\frac{i+1}{2}]}\binom{q/2-1}{[i/2]}\right).$$

The reader may verify that $A_1 = A_2 = 0$ and that $A_3 = (q-1)(q-2)/6$, showing once more that the minimum distance of $\mathfrak{H}_m$ is 3. For large $q$, the main contribution to $A_i$ is $\frac{1}{q}\binom{q-1}{i}$. The other term is much smaller.

## 4. THE DOUBLE ERROR CORRECTING BCH CODES

Let $q = 2^m > 4$. We define the double error correcting BCH codes $\mathfrak{B}_m$ to be the ideals in $\mathbf{F}_2[T]/(T^{q-1}-1)$ generated by $f_\alpha(T)f_{\alpha^3}(T)$ where $\alpha$ is a generator of $\mathbf{F}_q^*$. The codes $\mathfrak{B}_m$ are very classical. Their definition does not depend on the choice of $\alpha$. The dimension of the codes $\mathfrak{B}_m$ is $q - 1 - 2m$ and their minimum distance is 5. See [8, p. 201] for a proof. The smallest non-trivial example is the ideal generated by the polynomial

$$T^8 + T^7 + T^6 + T^4 + 1 = (T^4 + T + 1)(T^4 + T^3 + T^2 + T + 1)$$

in the ring $\mathbf{F}_2[T]/(T^{15} - 1)$. This code has length 15 and dimension 7.

To determine the weight distributions of $\mathfrak{B}_m$, we proceed as we did with the Hamming codes $\mathfrak{H}_m$ in the previous section. First, since $\alpha$ generates the multiplicative group, $\alpha^3$ is not a Galois conjugate of $\alpha$. Therefore the code $\mathfrak{B}_m$ is the restriction of the ideal $(T - \alpha)(T - \alpha^3) \subset \mathbf{F}_q[T]/(T^{q-1}-1)$

$$\begin{array}{ccc} \mathfrak{B}_m & \xleftarrow{\text{Res}} & (T-\alpha)(T-\alpha^3) \\ \Big\uparrow{\text{dual}} & & \Big\uparrow{\text{dual}} \\ \mathfrak{B}_m^\perp & \xleftarrow{\text{Trace}} & (T-\alpha)(T-\alpha^3)^\perp \end{array}$$

The polynomials in $(T - \alpha)(T - \alpha^3)$ are precisely those that vanish in $\alpha$ and $\alpha^3$. Arguing as in the previous section, we conclude that $(T - \alpha)(T - \alpha^3)^\perp$ is given by

$$\{(\lambda x + \mu x^3)_{x \in \mathbf{F}_q^*} : \lambda, \mu \in \mathbf{F}_q\}$$

and therefore, by Delsarte's Theorem,

$$\mathfrak{B}_m^\perp = \{(\text{Tr}(\lambda x + \mu x^3))_{x \in \mathbf{F}_q^*} : \lambda, \mu \in \mathbf{F}_q\}.$$

We want to compute the weight distributions of the codes $\mathfrak{B}_m^\perp$. We proceed as in [14]. Let $\mathbf{v}_{\lambda,\mu} = (\text{Tr}(\lambda x + \mu x^3))_{x \in \mathbf{F}_q^*}$ be a word in $\mathfrak{B}^\perp$. Its weight is closely related to the number of rational points on the curve

$$E: \quad Y^2 - Y = \lambda X + \mu X^3.$$

Over $\mathbf{F}_q$: for $x \in \mathbf{F}_q$ one has that $\text{Tr}(\lambda x + \mu x^3) = 0$ if and only if there exists $y \in \mathbf{F}_q$ with $y^2 - y = \lambda x + \mu x^3$. Therefore, taking care of the points $(0,0)$, $(0,1)$, and the point at infinity, we have that

$$|\mathbf{v}_{\lambda,\mu}| = (q-1) - \frac{1}{2}(\#E(\mathbf{F}_q) - 3).$$

In retrospect we could say that, in a similar way, the Hamming codes of the previous section are related to the family of curves $Y^2 - Y = \nu X$, $\nu \in \mathbf{F}_q$. Since these are all curves of genus 0, the calculation in section 3 was particularly simple.

For $\lambda = \mu = 0$ one obtains the zero word. If $\mu = 0$ and $\lambda \neq 0$, the curve $E$ is a conic; this gives rise to $q - 1$ words of weight $q - 1 - \frac{1}{2}(q + 1 - 3) = q/2$.

When $\mu \neq 0$, the curve $E$ has genus 1. It is easily seen to be a supersingular elliptic curve. This means that it has no points of order 2 or, equivalently, that its ring of $\overline{\mathbf{F}}_2$-endomorphisms is not commutative. See Deurings beautiful 1941 paper [2] for a discussion of supersingular elliptic curves. There is, up to isomorphism, only one supersingular elliptic curve over $\overline{\mathbf{F}}_2$. It has equation $Y^2+Y = X^3$. Over $\mathbf{F}_q$ there are more isomorphism classes. They are all isomorphic over $\overline{\mathbf{F}}_2$. The following proposition describes the situation.

**Proposition 4.1.** *Let $\mathbf{F}_q$ be an extension of $\mathbf{F}_2$ of degree $m$. If $m$ is odd, there are, up to isomorphism, exactly 3 supersingular curves over $\mathbf{F}_q$. If $m$ is even, there are, up to isomorphism, exactly 7 supersingular curves over $\mathbf{F}_q$. In the following tables they are sorted according to their number of points over $\mathbf{F}_q$. We also give the cardinality of their groups of $\mathbf{F}_q$-rational automorphisms:*

$m$ odd

| $\#E(\mathbf{F}_q)$ | freq. | $\#\mathrm{Aut}_{\mathbf{F}_q}(E)$ |
|---|---|---|
| $q+1-\sqrt{2q}$ | 1 | 4 |
| $q+1$ | 1 | 2 |
| $q+1+\sqrt{2q}$ | 1 | 4 |

$m$ even

| $\#E(\mathbf{F}_q)$ | freq. | $\#\mathrm{Aut}_{\mathbf{F}_q}(E)$ |
|---|---|---|
| $q+1-2\sqrt{q}$ | 1 | 24 |
| $q+1-\sqrt{q}$ | 2 | 6 |
| $q+1$ | 1 | 4 |
| $q+1+\sqrt{q}$ | 2 | 6 |
| $q+1+2\sqrt{q}$ | 1 | 24 |

*Proof.* See [2, 15, 11]. It is actually very easy to give explicit equations for the curves. For instance, when $m$ is odd, the three curves are given by the equations $Y^2+Y = X^3$, $Y^2+Y = X^3+X$, and $Y^2+Y = X^3+X+1$ with $q+1$, $q+1\pm\sqrt{2q}$, and $q+1\mp\sqrt{2q}$ rational points respectively. Here $\pm$ denotes $(-1)^{(m^2-1)/8}$. We do not need these equations.

Now we count how often an isomorphism class occurs in our set of curves:

**Proposition 4.2.** *Let $E$ be a supersingular elliptic curve over $\mathbf{F}_q$. The number of curves in the family*
$$Y^2 + Y = \lambda X + \mu X^3 \qquad (\lambda \in \mathbf{F}_q,\, \mu \in \mathbf{F}_q^*)$$
*that are isomorphic over $\mathbf{F}_q$ to $E$, is equal to*
$$\frac{(q-1)(\#E(\mathbf{F}_q)-1)}{\#\mathrm{Aut}_{\mathbf{F}_q}(E)}.$$

*Proof.* If in the equation we replace $X$ by $\mu^{-1}X$ and $Y$ by $\mu^{-1}Y$ and then substitute $\lambda$ by $\lambda\mu^{-1}$, our family of supersingular elliptic curves becomes
$$Y^2 + \mu Y = X^3 + \lambda X \qquad (\lambda \in \mathbf{F}_q,\, \mu \in \mathbf{F}_q^*).$$
Let us fix one such curve $E$. A transformation of the form
$$Y \leftarrow u^3 Y + rX + s,$$
$$X \leftarrow u^2 X + t,$$
with $r,s,t \in \mathbf{F}_q$, $u \in \mathbf{F}_q^*$, maps $E$ to another curve in the family if and only if
$$s^2 + \mu s = t^3 + \lambda t,$$
$$t = u^{-4} r^2.$$
There are $(\#E(\mathbf{F}_q)-1)(q-1)$ such transformations, and all the corresponding curves are isomorphic to $E$. Since *all* $\mathbf{F}_q$-automorphisms of $E$ are also of this form, the number of curves in the family isomorphic to $E$ is
$$\frac{(q-1)(\#E(\mathbf{F}_q)-1)}{\#\mathrm{Aut}_{\mathbf{F}_q}(E)}.$$

Summing this over the isomorphism classes of $E$, we find, using Proposition 4.1, precisely $q(q-1)$ curves, i.e. we obtain all curves in the family. This proves the proposition.

From Propositions 4.1 and 4.2 it is easy to obtain the weight distribution of the codes $\mathfrak{B}_m^\perp$:

$m$ odd

| weight | frequency |
|---|---|
| 0 | 1 |
| $\frac{q+\sqrt{2q}}{2}$ | $\frac{q-1}{4}(q-\sqrt{2q})$ |
| $\frac{q}{2}$ | $\frac{q-1}{2}q+(q-1)$ |
| $\frac{q-\sqrt{2q}}{2}$ | $\frac{q-1}{4}(q+\sqrt{2q})$ |

$m$ even

| weight | frequency |
|---|---|
| 0 | 1 |
| $\frac{q+2\sqrt{q}}{2}$ | $\frac{q-1}{24}(q-2\sqrt{q})$ |
| $\frac{q+\sqrt{q}}{2}$ | $\frac{q-1}{3}(q-\sqrt{q})$ |
| $\frac{q}{2}$ | $\frac{q-1}{4}q+(q-1)$ |
| $\frac{q-\sqrt{q}}{2}$ | $\frac{q-1}{3}(q+\sqrt{q})$ |
| $\frac{q-2\sqrt{q}}{2}$ | $\frac{q-1}{24}(q+2\sqrt{q})$ |

These distributions were already found in a different way by Kasami [5, 6]; see also [8, Ch. 8, sect. 7 and Ch. 15, sect. 4]. Applying the MacWilliams Identities, one obtains explicit formulas for the weight distributions of the BCH codes $\mathfrak{B}_m$ themselves. We leave this calculation to the reader.

## 5. The Melas codes

In 1960 the Melas codes [9] were introduced by C. M. Melas, who worked for IBM. The definition of these codes is very similar to the one of the BCH codes. Let $\mathbf{F}_q$ be an extension of $\mathbf{F}_2$. Assume $q = 2^m > 4$. The Melas codes are ideals in $\mathbf{F}_2[T]/(T^{q-1}-1)$ generated by $f_\alpha(T)f_{\alpha^{-1}}(T)$ where $\alpha$ is a generator of $\mathbf{F}_q^*$. The definition does not depend on the choice of $\alpha$. The dimension of the codes $\mathfrak{M}_m$ is $q-1-2m$ and their minimum distance is at least 3 or 5, depending on whether $m$ is even or not.

The smallest non-trivial example occurs for $m = 4$. It is the ideal generated by the 15th cyclotomic polynomial (mod 2)
$$\Phi_{15} = T^8 + T^7 + T^5 + T^4 + T^3 + T + 1 = (T^4 + T + 1)(T^4 + T^3 + 1)$$
in the ring $\mathbf{F}_2[T]/(T^{15}-1)$. The length of this code is 15 and its dimension is 7.

We proceed as we did with the BCH codes $\mathfrak{B}_m$ in the previous section. $\mathfrak{M}_m$ is the restriction of the ideal $(T-\alpha)(T-\alpha^{-1}) \subset \mathbf{F}_q[T]/(T^{q-1}-1)$. We conclude, as in the previous section, that the $\mathbf{F}_q$-dual of this code is
$$\{(\lambda x + \mu x^{-1})_{x \in \mathbf{F}_q^*} : \lambda, \mu \in \mathbf{F}_q\},$$
and therefore, by Delsarte's Theorem,
$$\mathfrak{M}_m^\perp = \{(\mathrm{Tr}(\lambda x + \mu x^{-1}))_{x \in \mathbf{F}_q^*} : \lambda, \mu \in \mathbf{F}_q\}.$$
Let $\mathbf{v}_{\lambda,\mu} = (\mathrm{Tr}(\lambda x + \mu x^{-1}))_{x \in \mathbf{F}_q^*}$ be a word in $\mathfrak{M}_m^\perp$. If $\lambda = \mu = 0$, we find the zero word. If only one of $\lambda$, $\mu$ is zero, the code word has weight $q/2$. There are $2q-1$ such words. In the remaining cases the weight is closely related to the number of rational points on the curve
$$E: \quad Y^2 - Y = \lambda X + \mu X^{-1}$$

over $\mathbf{F}_q$. Multiplying the equation by $\lambda^2 X^2$ and replacing $Y$ by $X^{-1}Y$ and $X$ by $\lambda^{-1}X$, we find the equation of an elliptic curve

$$E: \quad Y^2 + XY = X^3 + \nu^2 X,$$

where $\nu$ is uniquely determined by $\nu^2 = \mu\lambda \in \mathbf{F}_q^*$. Keeping in mind the points with $X = 0$ and the point at infinity, one easily computes the weight of the corresponding code word

$$|\mathbf{v}_{\lambda,\mu}| = (q-1) - \frac{1}{2}(\#E(\mathbf{F}_q) - 2) = \frac{q-1+t}{2}$$

where $t$ is defined by $\#E(\mathbf{F}_q) = q+1-t$. Let us also note that the point $(\nu, 0)$ on $E$ has order 4 in the group $E(\mathbf{F}_q)$.

To determine the weight distribution of the codes $\mathfrak{M}^\perp$, we must count how many curves $E$ there are in our family with a given number of points. We first compute that the $j$-invariant of $Y^2 + XY = X^3 + \nu^2 X$ is equal to $\nu^{-4}$. This shows that for every non-zero $j$-invariant there is exactly one elliptic curve in the family and this elliptic curve is not supersingular. Up to isomorphism, there are exactly two elliptic curves over $\mathbf{F}_q$ with $j$-invariant equal to $\nu^{-4}$. If one has $q+1-t$ points over $\mathbf{F}_q$, the other has $q+1+t$ such points. We can decide which curve occurs by observing that the point $(\nu, 0) \in E(\mathbf{F}_q)$ has order 4. Therefore, the curve with $q+1-t$ points, where $t \equiv q+1 \pmod 4$ occurs.

The following proposition describes how many elliptic curves there are with a given number of points over $\mathbf{F}_q$. We consider only non-supersingular elliptic curves or, equivalently, curves with an even number of points. Proposition 5.2 has been known since the work of Deuring [2]. It involves certain *class numbers.*

**Definition.** Let $d \in \mathbf{Z}$ be a negative integer congruent to 0 or 1 (mod 4). The *Kronecker class number* $H(d)$ is the cardinality of the set of $\mathrm{SL}_2(\mathbf{Z})$-orbits of positive definite binary quadratic forms $aX^2 + bXY + cY^2$ with $a, b, c \in \mathbf{Z}$ with discriminant $b^2 - 4ac$ equal to $d$.

It is well known that

$$H(d) = \#\{(a, b, c) \in \mathbf{Z}^3 : b^2 - 4ac = d, |b| \le a \le c \text{ and } b > 0$$
$$\text{whenever } |b| = a \text{ or } a = c\}.$$

For a triple $(a, b, c)$ that satisfies the conditions in the formula for $H(d)$ we have that $|d| = 4ac - b^2 \ge 4a^2 - a^2 = 3a^2$ and hence $|b| \le a \le \sqrt{|d|/3}$. This shows that $H(d)$ is finite. When $d$ is not too large, the formula is suitable for calculating $H(d)$ in practice.

TABLE 5.1

| $-d$ | $H(d)$ | $-d$ | $H(d)$ | $-d$ | $H(d)$ | $-d$ | $H(d)$ | $-d$ | $H(d)$ |
|---|---|---|---|---|---|---|---|---|---|
| 3 | 1 | 23 | 3 | 43 | 1 | 63 | 5 | 83 | 3 |
| 4 | 1 | 24 | 2 | 44 | 4 | 64 | 4 | 84 | 4 |
| 7 | 1 | 27 | 2 | 47 | 5 | 67 | 1 | 87 | 6 |
| 8 | 1 | 28 | 2 | 48 | 4 | 68 | 4 | 88 | 2 |
| 11 | 1 | 31 | 3 | 51 | 2 | 71 | 7 | 91 | 2 |
| 12 | 2 | 32 | 3 | 52 | 2 | 72 | 3 | 92 | 6 |
| 15 | 2 | 35 | 2 | 55 | 4 | 75 | 3 | 95 | 8 |
| 16 | 2 | 36 | 3 | 56 | 4 | 76 | 4 | 96 | 6 |
| 19 | 1 | 39 | 4 | 59 | 3 | 79 | 5 | 99 | 3 |
| 20 | 2 | 40 | 2 | 60 | 4 | 80 | 6 | 100 | 3 |

**Proposition 5.2.** *Let $\mathbf{F}_q$ be an extension of $\mathbf{F}_2$. Let $t \in \mathbf{Z}$ be an odd integer. The number $M(t)$ of elliptic curves over $\mathbf{F}_q$ up to isomorphism with $q+1-t$ points is given by*
$$M(t) = \begin{cases} H(t^2 - 4q), & \text{if } |t| < 2\sqrt{q}, \\ 0 & \text{otherwise.} \end{cases}$$

*Proof.* This is well known. One can find a proof in [2] or [15].

As a consequence we find the following weight distributions for the duals of the Melas codes [7]: the number of code words in $\mathfrak{M}_m^\perp$ of weight $(q+1-t)/2$ is given by
$$\begin{cases} (q-1)H(t^2 - 4q), & \text{if } t \in \mathbf{Z} - \{1\}, t \equiv q+1 \pmod{4} \text{ and } t^2 < 4q; \\ (q-1)(H(1-4q) + 2), & \text{if } t = 1; \\ 1, & \text{if } t = q+1 \text{ (the zero word)}; \\ 0, & \text{otherwise.} \end{cases}$$

To obtain the weight distributions of the Melas codes themselves, we apply Theorem 2.2. The MacWilliams Identities give complicated formulas for the weights of the Melas codes, involving sums of the type $\sum_t t^i H(t^2 - 4q)$. Luckily the Eichler-Selberg formula for the traces of the Hecke operator $T_q$ acting on cusp forms for the congruence subgroup
$$\Gamma_1(4) = \left\{ \begin{pmatrix} a & b \\ c & d \end{pmatrix} \in \mathrm{SL}_2(\mathbf{Z}) : a \equiv 1 \pmod{4} \text{ and } c \equiv 0 \pmod{4} \right\}$$
involves similar sums.

**Proposition 5.3** (Eichler-Selberg Trace Formula for $\Gamma_1(4)$). *Let $m \geq 1$ and $q = 2^m$. The trace of the Hecke operator $T_q$ acting on the space of cusp forms $S_k(\Gamma_1(4))$ of weight $k \geq 2$ is given by*
$$\mathrm{Tr}(T_q) = -1 + q - \sum_t H(t^2 - 4q) = 0 \qquad \text{when } k = 2,$$
$$= -1 - (-1)^{kq/2} \sum_t Q_{k-2}(t, q) H(t^2 - 4q) \qquad \text{otherwise.}$$

*Here $t$ runs over $\{t \in \mathbf{Z} : t^2 < 4q, \text{ and } t \equiv q+1 \pmod{4}\}$. The $Q_\kappa(t, n)$ are recursively defined by*
$$Q_0(t, n) = 1,$$
$$Q_1(t, n) = t,$$
$$Q_{\kappa+1}(t, n) = tQ_\kappa(t, n) - nQ_{\kappa-1}(t, n) \quad \text{for } \kappa \geq 1.$$

*One has that $Q_\kappa(t, n) = (\rho^{\kappa+1} - \overline{\rho}^{\kappa+1})/(\rho - \overline{\rho})$ where $\rho, \overline{\rho}$ denote the zeroes of the polynomial $X^2 - tX + n$.*

*Proof.* See [10] or [3].

The group $\Gamma_1(4)$ is closely related to the modular curve $X_1(4)$ which parametrizes elliptic curves together with a point of order 4. Since all curves in our family of elliptic curves admit a rational point of order 4, one would expect to find an approach which is more direct than our method using the MacWilliams Identities and the Trace formula.

Combining Proposition 5.3 with the MacWilliams Identities gives, after some elementary and tedious calculations, the following formula for the number $A_i$ of words in the Melas code $\mathfrak{M}_m$ of weight $i$:

**Theorem 5.4.**
$$q^2 A_i = \binom{q-1}{i} + 2(-1)^{[\frac{i+1}{2}]}(q-1)\binom{q/2-1}{[i/2]}$$
$$- (q-1) \sum_{\substack{j=0 \\ j \equiv i \pmod{2}}}^{i} W_{i,j}(q)(1+\tau_{j+2}(q));$$

here the polynomials $W_{i,j}(q)$ are zero whenever $j<0$ or $j>i$. For $0 \leq j \leq i$ and $i \equiv j \pmod{2}$ they are defined by

$$W_{0,0} = 1, \qquad W_{1,1} = -1$$
$$(i+1)W_{i+1,j+1} = -qW_{i,j+2} - W_{i,j} - (q-i)W_{i-1,j+1};$$

by $\tau_\kappa(q)$ we denote for $\kappa \geq 3$ the trace of the Hecke operator $T_q$ on the space of cusp forms $S_\kappa(\Gamma_1(4))$. For convenience we let $\tau_2(q) = -q$.

Using Theorem 5.4, it is not difficult to calculate the frequencies of the small weights. This has been done in [12]. The computation naturally splits into two parts: for even weights (in the sense of coding theory), one only encounters traces of Hecke operators acting on cusp forms of even weight (in the sense of modular forms). This involves only modular forms for the group $\Gamma_0(4)$. By means of the theory of Atkin and Lehner [1], one reduces the calculation to the computation of the traces of $T_q$ for the full modular group $SL_2(\mathbf{Z})$. Similarly, for odd weights one encounters only cusp forms of odd weight. In this case cusp forms for the group $\Gamma_0(4)$ with respect to the unique non-trivial character of conductor 4 are involved. In this case all forms are "new" in the sense of Atkin and Lehner. We refer to [12] for the details and several explicit calculations.

## References

1. Atkin, A.O.L. and Lehner, J., *Hecke operators on* $\Gamma_0(m)$, Math. Ann. **185** (1970), 134–160.
2. Deuring, M., *Die Typen der Multiplikatorenringe elliptischer Funktionenkörper*, Abh. Math. Sem. Hansischen Univ. **14** (1941), 197–272.
3. Cohen, H., *Trace des opérateurs de Hecke sur* $\Gamma_0(N)$, Sém. Théorie Nombres **4** (1976–1977), Bordeaux.
4. Goppa, V.G., *Codes on algebraic curves*, Soviet Math. Dokl. **24** (1981), 170–172.
5. Kasami, T., *Weight distributions of Bose-Chaudhuri-Hocquenghem codes* (Bose, R.S. and Dowling, T.A., eds.), Combinatorial Math. and its Applications, Univ. of North Carolina Press, Chapel Hill, NC, 1969.
6. ______, *The weight enumerators for several classes of subcodes of the 2nd order binary Reed-Muller codes*, Inform. and Control (Shenyang) **18**, 369–394.
7. Lachaud, G. and Wolfmann, J., *The weights of the orthogonals of the extended quadratic binary Goppa codes*, IEEE Trans. Inform. Theory **36** (1990), 686–692.
8. MacWilliams, J. and Sloane, N.J.A., *The theory of error-correcting codes*, North Holland, Amsterdam, 1983.
9. Melas, C.M., *A cyclic code for double error correction*, IBM J. Res. Develop. **4** (1960), 364–366.
10. Oesterlé, J., *Sur la trace des opérateurs de Hecke*, Thèse de 3° cycle, Orsay, 1977.
11. Schoof, R., *Non-singular plane cubic curves over finite fields*, J. Combin. Theory Ser. A **46** (1987), 183–211.

12. Schoof, R. and Van der Vlugt, M., *Hecke operators and the weight distributions of certain codes*, J. Combin. Theory Ser. A **57** (1991), 163–186.
13. Silverman, J., *The arithmetic of elliptic curves*, Graduate Texts in Math., vol. 106, Springer-Verlag, New York, 1986.
14. Van der Geer, G. and Van der Vlugt, M., *Reed-Muller codes and supersingular curves*. I, Compositio Math. **84** (1992), 333–367.
15. Waterhouse, W., *Abelian varieties over finite fields*, Ann. Sci. École Norm. Sup. **2** (1969), 521–560.

DIPARTIMENTO DI MATEMATICA, 2A UNIVERSITÀ DI ROMA "TOR VERGATA", I-00133 ROMA, ITALY
*E-mail address*: schoof@volterra.science.unitn.it